\documentclass{amsart}
\usepackage{amsrefs}
\newtheorem{theorem}{Theorem}
\newtheorem{cor}{Corollary}

\title{A very short proof of Cauchy's interlace theorem for eigenvalues of Hermitian matrices}
\author{Steve Fisk \\ Bowdoin College \\ fisk@bowdoin.edu}

\begin{document}
\maketitle

We use an overlooked characterization of interlacing to give a two
sentence proof of Cauchy's interlace theorem\cite{hwang}.  Recall that
if polynomials $f(x)$ and $g(x)$ have all real roots $r_1\le
r_2\le\dots\le r_n$ and $s_1\le s_2\le\dots\le s_{n-1}$ then we say
that $f$ and $g$ \emph{interlace} if and only if
$$ r_1 \le s_1 \le r_2 \le s_2 \cdots \le s_{n-1} \le r_n.$$
The following can be
found in \cite{rahman}, along with a discussion of its history back to
Hermite.

\begin{theorem}
  The roots of polynomials $f,g$ interlace if and only if  the
  linear combinations $f+\alpha g$ have all real roots for all
  $\alpha\in{{\mathbb R}}$.
\end{theorem}

\begin{cor}
  If $A$ is a Hermitian matrix, and $B$ is a principle submatrix of
  $A$, then the eigenvalues of $B$ interlace the eigenvalues of $A$.
\end{cor}

\begin{proof}
Choose $\alpha\in\mathbb R$,  partition $
A =\left(\begin{array}{c|c}
B & c \\\hline c^\ast & d
\end{array}\right)
$
and consider the following equation that follows from linearity of the
determinant:

$$
\left|\begin{array}{c|c}
  B - xI  & c \\\hline c^\ast & d-x+\alpha
\end{array}\right|
=
\left|\begin{array}{c|c}
  B - xI  & c \\\hline c^\ast & d-x
\end{array}\right|
+
\left|\begin{array}{c|c}
  B - xI  & c \\\hline 0 & \alpha
\end{array}\right|
$$

Since the matrix on the left hand side is the characteristic
polynomial of a Hermitian matrix, $|A-xI| + \alpha|B-xI|$ has all real
roots for any $\alpha$, and hence the eigenvalues interlace.
\end{proof}

\end{document}